\documentclass{article}

\usepackage[includeheadfoot,
            bindingoffset=0mm,
            inner   = 25mm,
            top     = 20mm,
            outer   = 25mm,
            bottom  = 10mm,
            paperwidth = 210mm,
            paperheight = 297mm,
            ]{geometry}
\usepackage[margin={2.0cm,0cm},oneside,labelfont={sf,bf},singlelinecheck=false]{caption}

\usepackage{graphicx}

\usepackage[fleqn]{amsmath}
\usepackage{accents}
\usepackage{amssymb}
\usepackage{amsmath}

\usepackage{enumitem}
\usepackage{cite}
\usepackage[perpage]{footmisc}
\usepackage{hyphenat}
\usepackage[english]{babel}
\usepackage[section]{placeins}
\usepackage{color}
\usepackage{upgreek}
\usepackage{listings}

\definecolor{codegreen}{rgb}{0,0.6,0}
\definecolor{codeblue}{rgb}{0,0,0.8}
\definecolor{codegrey}{rgb}{0.5,0.5,0.5}
\usepackage{amsthm}
\theoremstyle{definition}

\newtheorem{remark}{Remark}
\newtheorem{algorithm}{Algorithm}

\usepackage{titlesec}
\titleformat{\section}[hang]{\Large\bfseries\raggedright\sffamily}{\thesection}{1em}{}
\titleformat{\subsection}[hang]{\large\bfseries\raggedright\sffamily}{\thesubsection}{1em}{}
\titleformat{\subsubsection}[hang]{\normalsize\bfseries\raggedright\sffamily}{\thesubsubsection}{1em}{}
\usepackage{abstract}

\usepackage{authblk}

\begin{document}

\title{ \huge\bfseries\sffamily Canonical block-oriented model }

\author[1]{A. Polar}
\author[2]{M. Poluektov}
\affil[1]{Independent Software Consultant, Duluth, GA, USA}
\affil[2]{International Institute for Nanocomposites Manufacturing, WMG, University of Warwick, Coventry CV4 7AL, UK}

\date{ \huge\normalfont\sffamily DRAFT: \today }

\maketitle

\setlength{\absleftindent}{2.0cm}
\setlength{\absrightindent}{2.0cm}
\setlength{\absparindent}{0em}
\begin{abstract}
The block-oriented models are usually based on linear dynamic and non-linear static blocks that are connected in various sequential/parallel ways. Some particular configurations of the involved blocks result in the well-known Hammerstein, Wiener, Hammerstein-Wiener and generalised Hammerstein models. The Urysohn model is a lesser-known model; it is represented by a single non-linear dynamic block and can be approximated by a number of parallel Hammerstein blocks. In this paper, it is shown that any block-oriented model can be adequately replaced by a single Urysohn block followed by a single static non-linear block. Furthermore, a method of the so-called non-parametric identification of such object is introduced. \\
\textbf{Keywords:} Urysohn model, non-linear system identification, block-oriented models, non-parametric identification.
\end{abstract}

\section{Introduction}
\label{sec:intro}

Modelling dynamic objects using sequential or parallel connection of linear dynamic and non-linear static blocks is very common. Such representation has sufficient descriptive abilities for most applications, while maintaining a simple structure. The classical models of this type are the Hammerstein, the Wiener, the Hammerstein-Wiener and the generalised Hammerstein\footnote{The generalised Hammerstein model is represented by a finite number of parallel Hammerstein blocks.} models. In addition to the aforementioned models, there is also the Urysohn model, which is represented by a single non-linear dynamic block. 

Having a discrete-time Wiener model as 
\begin{align}
  &y_i = \sum_{j=1}^{m} {h_j x_{i-j+1}} ,
  \label{eq:Wiener1} \\
  &z_i = f\left(y_i\right) ,
  \label{eq:Wiener2}
\end{align}
where $x_k$ is the input, $y_k$ and $z_k$ are the intermediate variable and the final output, respectively, $f$ is the static nonlinearity, $h_j$ are the linear model coefficients, it is possible to generalise it to a more complex object by an elementary modification. At the first step, each linear term in the linear block is replaced by a non-linear term with functions $g_j$, 
\begin{equation}
  y_i = \sum_{j=1}^{m} {g_j\left(x_{i-j+1}\right)} .
  \label{eq:Ur2}
\end{equation}
At the second step, static nonlinearity $f$ is converted into a sum of nonlinearities by including a few preceding intermediate variables,
\begin{equation}
  z_i = \sum_{j=1}^{p} {f_j\left(y_{i-j+1}\right)} .
  \label{eq:Ur1}
\end{equation}
At both steps, the one can see the structure of the Urysohn model, as both equations \eqref{eq:Ur2} and \eqref{eq:Ur1} are non-linear dynamic blocks of a discrete-time system. Thus, the obtained model \eqref{eq:Ur2}-\eqref{eq:Ur1} becomes the system of two sequential discrete-time Urysohn operators, or the two-Urysohn model for simplicity.

It is easy to show that the two-Urysohn model cannot be replaced accurately by a single Urysohn block. When a single Urysohn block is considered, e.g. \eqref{eq:Ur2}, and functions $g_j$ are expressed by polynomials, the terms that contribute to output $y_i$ can be represented as powers of input values $x_k^q$. However, when two sequential Urysohn blocks are considered, the intermediate variable, which already contains sums of powers of inputs $x_k^q$, is converted by the second Urysohn block. Therefore, terms that contribute to output $z_i$ will also contain products of time-shifted input values, such as $x_k^q x_l^r$, where $k \neq l$, which is the principal difference between the single Urysohn block and the sequence of such blocks.

The two-Urysohn model is a generalisation of popular block-oriented models --- the Hammerstein, the Wiener, the Hammerstein-Wiener and the generalised Hammerstein models\footnote{The Hammerstein block can be written as 
\begin{equation*}
  y_i = \sum_{j=1}^{m} {h_j u\left(x_{i-j+1}\right)} ,
\end{equation*}
where $h_j$ are scalar parameters of a linear block and $u$ is some non-linear function. By defining $g_j\left(x\right) = h_j u\left(x\right)$, the Urysohn block \eqref{eq:Ur2} is obtained. Similarly, it can be shown that the generalised Hammerstein is also a particular case of the Urysohn model.}. There is a wide variety of identification methods of the aforementioned popular models and there are comprehensive review papers on this topic \cite{Schoukens2019}; however, the two-Urysohn model has never been identified before. The aim of this paper is to analyse the structure of the two-Urysohn model, to exclude possible redundancy from this model and to suggest an identification algorithm for the model.

Identification algorithm for the single Urysohn model has already been proposed by the authors of this paper in \cite{Poluektov2019}, where properties of the single Urysohn model have been analysed as well. To facilitate understanding of the theory of this paper, major concepts from \cite{Poluektov2019} are repeated in sections \ref{sec:single} and \ref{sec:basic}. This is followed by the analysis of the two-Urysohn model, section \ref{sec:canonical}, and the identification algorithm for the two-Urysohn model, section \ref{sec:canonicalidentification}. Finally, computational examples are given in section \ref{sec:simulation}.

\section{Forms of the discrete-time Urysohn model}
\label{sec:single}

The discrete-time Urysohn model, e.g. \eqref{eq:Ur2}, contains non-linear functions $g_j$ that must be represented in some form, before the model can be identified or used to reproduce the input-output relation of a control system. These functions can be represented by polynomials with some coefficients \cite{Makarov1994,Makarov2012} or even by a set of arbitrary functions with some weights \cite{Schoukens2011}, which results in a representation of the Urysohn model by parallel Hammerstein blocks.

The authors of this paper have previously proposed to take functions $g_j$ to be either piecewise constant or piecewise linear \cite{Poluektov2019}, which gives two possible representations. A set of functions $g_j$ can be called a kernel of the discrete-time Urysohn model, similarly to the kernel of the continuous-time model, which is given by an integral equation. For the piecewise constant and linear kernels, the input range $\left[ x_\mathrm{min}, x_\mathrm{max} \right]$ is divided into a number of equal intervals, within which functions $g_j$ are either constant or linear. The advantages of such representations of $g_j$ are the controllable accuracy of the model (by varying the density of grid points) and the simplicity of the identification of the model.

From the mathematical point of view, it is convenient to collect all grid point values of $g_j$ in a rectangular matrix and rewrite model \eqref{eq:Ur2} in the corresponding way. The piecewise-constant-kernel (PCK) discrete-time Urysohn model becomes \cite{Poluektov2019}
\begin{align}
  &y_i = \sum_{j=1}^{m} U\left[ j, k_{i-j+1} \right] ,
  \label{eq:UryshonDiscr} \\
  &k_i = 1 + \operatorname{round}\left( \frac{ \left(n-1\right) \left(x_i-x_\mathrm{min}\right) }
    { x_\mathrm{max} - x_\mathrm{min} } \right) ,
  \label{eq:UryshonDiscrContr}
\end{align}
where $U$ is the matrix with indices shown in $\left[ \cdot, \cdot \right]$, operator $\operatorname{round}\left(\cdot\right)$ is the rounding to the nearest integer and $n$ is the number of grid points for $g_j$. 

The case of piecewise constant $g_j$ can also be referred to as the quantised Urysohn model, as it is equivalent to quantisation of the input first and subsequent calculation of the output. The quantised discrete-time Urysohn operator was originally investigated by one of the authors of this paper in 1990s \cite{Poluektov1990} and successfully applied to a real-time modelling of diesel engines. Since then, the authors did not find any similar approaches.

The case of piecewise linear $g_j$ can be rewritten as \cite{Poluektov2019}
\begin{align}
  &y_i = \sum_{j=1}^{m} \left( \left( 1 - \psi_{i-j+1} \right) U\left[ j, k^\mathrm{o}_{i-j+1} \right] +
    \psi_{i-j+1} U\left[ j, k^\mathrm{v}_{i-j+1} \right] \right) , \label{eq:UryshonLin} \\
  &k_i^\mathrm{o} = \left\lfloor b_i \right\rfloor , \quad
    k_i^\mathrm{v} = \left\lceil b_i \right\rceil , \quad
    \psi_i = b_i - k_i^\mathrm{o} , \\
  &b_i = 1 + \left( n-1 \right) \frac{x_i - x_\mathrm{min}}{x_\mathrm{max} - x_\mathrm{min}} , \label{eq:UryshonLinContr}
\end{align}
where $\lfloor \cdot \rfloor$ and $\lceil \cdot \rceil$ are the floor and the ceiling functions, respectively. It is useful to refer to this model as the piecewise-linear-kernel (PLK) discrete-time Urysohn model. In both PCK and PLK cases, models are represented on grids and the only difference is the computation of the output.

\section{Identification of the single Urysohn block}
\label{sec:basic}

The idea of the identification method for the single Urysohn block relies on real-time modification of matrix $U$, such that the model output becomes equal to the measured output. According to models presented in section \ref{sec:single}, each input fragment of length $m$ defines a number of elements of matrix $U$, a weighted sum of which must be equal to the corresponding output value. In the case when this sum is not equal to the output (i.e. there is a discrepancy between the model and the measurement), each involved element of the matrix is corrected to fit this particular output. By repeating this numerical update of matrix $U$ for long enough input/output sequence, the model is tuned to reproduce the actual output. As shown in \cite{Poluektov2019}, this algorithim is a particular case of the projection descent method \cite{Kaczmarz1937,Tewarson1969} with quick convergence and good error filtering capabilities. The full algorithm is given below for the convenience of the reader.

\begin{algorithm}\label{alg:algOne}
Algorithm for identification of the PCK and PLK discrete-time Urysohn operators --- reproduced from \cite{Poluektov2019}.
\begin{enumerate}
  \item Assume initial approximation of matrix $U$, which can be arbitrary, including the all-zero matrix.
  \item Start with $i = m$.
  \item Calculate model output $\hat{y}_i$ based on actual inputs $\left( x_i, x_{i-1}, \ldots, x_{i-m+1} \right)$ and the current approximation of matrix $U$. \label{item:iterS}
  \item Calculate difference $D = y_i - \hat{y}_i$, where $y_i$ is the actual recorded output and $\hat{y}_i$ is the model output.
  \item Modify matrix $U$, such that a value proportional to $D$ is added to each element that was involved in the calculation of $\hat{y}_i$. \label{item:iterE}
  \item Increase index $i$ by $1$ and repeat steps \ref{item:iterS}-\ref{item:iterE} until $D$ becomes sufficiently small for sufficiently large number of iterations consecutively.
\end{enumerate}
\end{algorithm}

For the PCK model, the value referred to at step \ref{item:iterE} is $\alpha D / m$ and it is added to $\left[ j, k_{i-j+1} \right]$ elements of matrix $U$. Parameter $\alpha \in \left(0, 1\right]$ is introduced for suppressing the noise. For the PLK model, at step \ref{item:iterE} of the algorithm, $\alpha D \left( 1 - \psi_{i-j+1} \right) / \bar{\psi}_i$ is added to elements $\left[ j, k^\mathrm{o}_{i-j+1} \right]$ of matrix $U$ and  $\alpha D \psi_{i-j+1} / \bar{\psi}_i$ is added to elements $\left[ j, k^\mathrm{v}_{i-j+1} \right]$ of matrix $U$, where
\begin{equation}
  \bar{\psi}_i = \sum_{j=1}^{m} \left( \left( 1 - \psi_{i-j+1} \right)^2 + {\psi_{i-j+1}}^2 \right) .
\end{equation}

\section{Canonical block-oriented model}
\label{sec:canonical} 

The single Urysohn operator has certain descriptive limitations. Although it is a generalisation of the Hammerstein model and even an adequate replacement of multiple parallel Hammersteins, it cannot represent accurately Wiener objects, and the Wiener model, in turn, is not general enough to represent an Urysohn object. Here, the aim is to propose a block-oriented model, which can be used for every deterministic SISO (single input, single output) object with a finite memory, and which covers all known block-oriented models as particular cases.

As shown in the introduction, the model consisting of two sequential Urysohn blocks is a suitable candidate for such model. However, it contains a redundancy (discussed in section \ref{sec:descr}) and can be simplified further to a single Urysohn block followed by a static nonlinearity:
\begin{align}
  &y_i = \sum_{j=1}^{m} {g_j\left(x_{i-j+1}\right)} ,
  \label{eq:canonical1} \\
  &z_i = f\left(y_i\right) , \label{eq:canonical2}
\end{align}
Here, the Urysohn block, equation \eqref{eq:canonical1}, can be represented by either piecewise constant kernel, equations \eqref{eq:UryshonDiscr}-\eqref{eq:UryshonDiscrContr}, or piecewise linear kernel, equations \eqref{eq:UryshonLin}-\eqref{eq:UryshonLinContr}. The static nonlinearity, equation \eqref{eq:canonical2}, can be considered as a particular case of the discrete-time Urysohn model and can also be represented by either piecewise constant or linear kernels.

Formally, model \eqref{eq:canonical1}-\eqref{eq:canonical2} is capable of describing any deterministic SISO object with a finite memory with any predefined accuracy (shown in section \ref{sec:descr}); therefore, it can be called the canonical block-oriented model. However, in practice, it can occur that the structure of functions $g_j$ and $f$ is extremely complex and, thus, it is difficult to identify such model, i.e. functions $g_j$ and $f$ describing the object with a predefined accuracy may never be found, although the solution exists. This creates the case for the two-Urysohn model, equations \eqref{eq:Ur2}-\eqref{eq:Ur1}, redundancy of which allows for certain flexibility in tuning the kernels of the Urysohn blocks during the identification step and, thus, can lead to a more accurate identification. In the numerical examples of this paper, both the canonical and the two-Urysohn models are used.

\subsection{Descriptive capabilities}
\label{sec:descr}

Any deterministic SISO object with a finite memory can be represented in the discrete-time setting as 
\begin{equation}
  z_i = F\left( x_i, x_{i-1}, \ldots, x_{i-m+1} \right) ,
  \label{eq:SISO}
\end{equation}
where $F$ is some function of $m$ variables. The discussion in this section is limited to the case when $F$ is the continuous function of all variables\footnote{It is also possible to give the same proof for the case when $F$ contains a finite number of discontinuities of the first kind. This requires additional definition of the Urysohn model with piecewise-constant $g_j$, where the length of the intervals varies. For the purpose of this paper, the discussion of this case is omitted, as it overcomplicates the understanding of the presented concepts.}, which describes systems where small variation of the input results in a small variation of the output.

\textbf{Theorem.} The canonical model, equations \eqref{eq:canonical1}-\eqref{eq:canonical2}, is capable of describing any deterministic finite-memory SISO object with any predefined accuracy, given that output $z_i$ is a continuous function of all elements of the input sequence.

\textbf{Proof.} The proof consists of two parts. First, the statement of the theorem is proved for the quantised-input object. Second, the approximation of the continuous-input object by the quantised-input object is shown.

\textbf{I.} The quantised-input deterministic SISO object has input that takes integer values from $1$ to $n$. The object has the memory depth of $m$. Therefore, the exhaustive set of the quantised input sequences that determines all possible output values has the size of $n^m$. Assume such object is modelled using the canonical model with the special piecewise-constant-kernel Urysohn block. Kernel
\begin{equation}
  U = \begin{bmatrix}
    1       & 2        & 3        & \cdots & n      \\
    n       & 2n       & 3n       & \cdots & n^2    \\
    n^2     & 2n^2     & 3n^2     & \cdots & n^3    \\
    \vdots  & \vdots   & \vdots   & \ddots & \vdots \\
    n^{m-1} & 2n^{m-1} & 3n^{m-1} & \cdots & n^m    \\
  \end{bmatrix}
  \label{eq:Unm}
\end{equation}
provides a different output for every different input fragment from the exhaustive set. Sorted output values represent sequential integers from $\left( 1 + n + n^2 + \ldots + n^{m-1} \right)$ to $\left( n + n^2 + n^3 + \ldots + n^m \right)$. It is also useful to note that by subtracting $\left( n + n^2 + \ldots + n^{m-1} \right)$ from each possible output, the range limits are reduced to $\left[1, n^m\right]$, which means that there are exactly $n^m$ distinct outputs, the same number as the distinct input fragments.

Next, it is obviously possible to construct $f$ in \eqref{eq:canonical2} as a piecewise-constant function, which provides a unique mapping between every distinct value of intermediate variable $y_i$ and the corresponding final output $z_i$. For matrix $U$ given by \eqref{eq:Unm}, $y_i$ are integers and it is natural to call such piecewise-constant $f$ --- a lookup table. This means that the canonical model can map the exhaustive set of the quantised input fragments to the corresponding outputs and, therefore, can describe exactly any quantised-input object.

\textbf{II.} It is always possible to use a quantised-input object as an approximation for a continuous-input object, expressed by \eqref{eq:SISO} with $F$ being a continuous function, i.e. for any $\varepsilon > 0$ there exists $n$, such that $\left| z_i - \hat{z}_i \right| < 0$, where 
\begin{align*}
  &\hat{z}_i = F\left( \hat{x}_i, \hat{x}_{i-1}, \ldots, \hat{x}_{i-m+1} \right) , \\
  &\hat{x}_i = x_\mathrm{min} + \left( k_i - 1 \right) \frac{x_\mathrm{max} - x_\mathrm{min}}{n-1}
\end{align*}
and $k_i$ is given by equation \eqref{eq:UryshonDiscrContr}. Here, $\varepsilon$ is the predefined accuracy. \textbf{End of proof.}

The aim of this example is just to demonstrate this theoretical capability of the canonical model, even though using $U$ given by \eqref{eq:Unm} and $f$ in the form of a lookup table can be inconvenient from the practical point of view. It is also useful to note, that, in practice, when real physical objects are modelled, an approximation with a given input resolution is not an obstruction, since all signals are recorded in form of approximately-known instantaneous values and for high input resolution, the quantisation inaccuracy becomes insignificant compared to errors already introduced by sensor readings and the discrete time.

\subsection{Multiple inputs}

Although only SISO object has been considered so far, the descriptive capabilities of the Urysohn operator allow modelling objects with multiple inputs. In the case of the quantised kernel, it is possible to construct a mapping operator, which assigns a unique integer to any unique combination of input values\footnote{For example, in the case of two inputs and each input taking values $1$, $2$, $3$, a mapping matrix 
\begin{equation*}
  \begin{bmatrix}
    7 & 8 & 9  \\
    6 & 1 & 2  \\
    5 & 4 & 3  \\
  \end{bmatrix}
\end{equation*}
takes any pair of discrete inputs, which are used as row and column, and assigns a single unique value.}, which is plugged into model \eqref{eq:UryshonDiscr}. This means that there is no structural difference between a quantised-input SISO and a quantiesed-input MISO (multiple input, single output) Urysohn objects, the latter is simply larger in size. Further discussion of this property can be found in \cite{Poluektov2019}, where also a practical example of successful identification of a real physical object with two inputs is shown.

\section{Identification of the canonical block-oriented model}
\label{sec:canonicalidentification} 

The canonical model, equations \eqref{eq:canonical1}-\eqref{eq:canonical2}, is a particular case of the two-Urysohn model, equations \eqref{eq:Ur2}-\eqref{eq:Ur1}. Therefore, the algorithm of this section is presented for the two-Urysohn model, which also covers the canonical model.

Section \ref{sec:basic} defines formal steps to identify the single Urysohn block. The two-Urysohn model has two blocks and both must be identified using observable input $x_i$ and final output $z_i$. Intermediate variable $y_i$ is not observable; more than that, it is rather an auxiliary mathematical variable than a physical parameter. Obviously, in the case when $y_i$ is available, the identification is reduced to the already known case of the single Urysohn block. Indeed, having $y_i$ as the output of the first block and having $y_i$ as the input for the second block, both blocks can be identified using algorithm \ref{alg:algOne} of section \ref{sec:basic}. This idea is the foundation of the novel algorithm of identification of the two-Urysohn model. 

For provided input and output data and any approximation of the model in form of $g_j$ and $f_j$, intermediate value $y_i$ can be computed. Then, two new intermediate values are introduced --- incremented intermediate value $y_i + \Delta y$ and decremented intermediate value $y_i - \Delta y$. Having these three intermediate values, three output estimations are computed $\sum f_j\left(y_i\right)$, $\sum f_j\left(y_i + \Delta y\right)$, $\sum f_j\left(y_i - \Delta y\right)$. If either incremented or decremented value gives better approximation of final output $z_i$, it is accepted as the new intermediate value and both operators are corrected according to algorithm \ref{alg:algOne}. Applying this operation to every new pair of the observed input and output values, both operators converge to a stable form. 

Similar to a single Urysohn model, the output of the two-Urysohn model varies for each different initial approximation due to internal redundancy, which is a common property of many non-linear models, where identified parameters are included into products of sums.

Since the identification of the single Urysohn block works for both piecewise constant and linear forms, the blocks of the two-Urysohn model can be also in either of these forms. The formal algorithm is provided below.

\begin{algorithm}\label{alg:algTwo}
Algorithm for identification of the discrete-time two-Urysohn model.
\begin{enumerate}
  \item Assume initial approximation of functions $g_j$ and $f_j$, which can be arbitrary but not all zero functions. See remark \ref{rem:impr} below.
  \item Start with $i = m$.
  \item Calculate intermediate output $\hat{y}_i$ based on actual inputs $\left( x_i, x_{i-1}, \ldots, x_{i-m+1} \right)$ and the current approximation of $g_j$. \label{item:iter1}
  \item Introduce incremented and decremented intermediate values, $\hat{y}_i + \Delta y$ and $\hat{y}_i - \Delta y$. See remark \ref{rem:dlt} below. \label{item:iterI}
  \item Compute three possible model outputs for each intermediate value and select one of these three values, which provides the closest output to the observed output $z_i$.
  \item Use this selected intermediate value as an output for the first operator with functions $g_j$ and modify it according to steps \ref{item:iterS}-\ref{item:iterE} of algorithm \ref{alg:algOne}. See remark \ref{rem:realTime} below. \label{item:iter2m}
  \item Use this selected intermediate value as an input for the second operator with functions $f_j$ and modify it according to steps \ref{item:iterS}-\ref{item:iterE} of algorithm \ref{alg:algOne}. See remark \ref{rem:realTime} below. \label{item:iter2}
  \item Increase index $i$ by $1$ and repeat steps \ref{item:iter1}-\ref{item:iter2} until the difference between the computed final output and the actual final output becomes sufficiently small for sufficiently large number of iterations consecutively.
\end{enumerate}
\end{algorithm}

\begin{remark}\label{rem:impr}
In the presented form, the algorithm is the improvement of the existing initial approximation. Depending on the initial approximation, the algorithm can converge to different local minima, giving different errors, even for the exact input-output data. Therefore, from a practical point of view, it is reasonable to start from a number of initial guesses, run the algorithm in parallel, track their evolution and select a more accurate model at the end.
\end{remark}

\begin{remark}\label{rem:dlt}
The algorithm does not specify the exact value of $\Delta y$, optimal values of which might vary for different cases, as will be shown in numerical examples below. Furthermore, it might be varying from one iteration of the algorithm to another. The general idea is to keep $\Delta y$ small just to find the direction in which $\hat{y}_i$ must be changed. However, when more information regarding the optimal value of $\hat{y}_i$ is available, e.g. by knowing the structure of $f_i$, large $\Delta y$ might be taken.
\end{remark}

\begin{remark}\label{rem:realTime}
In the presented form, the algorithm operates in real-time and the operators are updated concurrently with data reading. In this case, it is sufficient to keep in memory only the $i$-th input and output values and the preceding $m-1$ input values. However, if long input and output sequences are already recorded, it is also possible to run the algorithm several times iterating through the data. In this case, it is possible to change only one of the operators during each run, i.e. execute either step \ref{item:iter2m} or step \ref{item:iter2} during one run and execute another step during a subsequent run. An example of such consecutive operator update will be discussed in section \ref{sec:simulation}.
\end{remark}

\section{Numerical examples}
\label{sec:simulation} 

The models introduced in this paper have internal redundancy. The variation in one part of the model can be compensated by another part. Therefore, accuracy of the identification method cannot be judged by comparison of the models, i.e. functions $g_j$ and $f_j$, and must be estimated by comparing the computed and the actual outputs. Usually, error measure
\begin{equation}
  E = \frac{\sqrt{\frac{1}{N} \sum_{i=1}^{N} \left(z_i - \hat{z_i}\right)^2}}{z_\mathrm{max} - z_\mathrm{min}} ,
  \label{eq:cost}
\end{equation}
where $z_i$ is actual output and $\hat{z_i}$ is modelled output, is used. Since the number of parameters is relatively large, there is a risk of overparametrisation. Therefore, error $E$ is computed for the validation (unseen) dataset. The input-output sequence used for the identification is usually called the training dataset. 

\subsection{Two sequential Urysohn blocks}
\label{sec:U2U}

The aim of this example is to show that the system consisting of two sequential Urysohn blocks can be identified using algorithm \ref{alg:algTwo}.

The idea of the numerical experiment is as follows. First, the original operators are generated. Next, the input sequence is generated. Afterwards, the output sequence is calculated using the original operators. This gives the precise (actual) input-output dependence, which is split into two parts --- the training and the validation datasets. Following this, the new operators are identified using algorithm \ref{alg:algTwo} and using just the training dataset. The obtained model is applied to the input of the validation dataset to obtain the modelled output. Finally, the accuracy of the obtained model is verified by calculating error $E$ using just the validation dataset.

\begin{figure}
  \begin{center}
    \includegraphics[width=6cm]{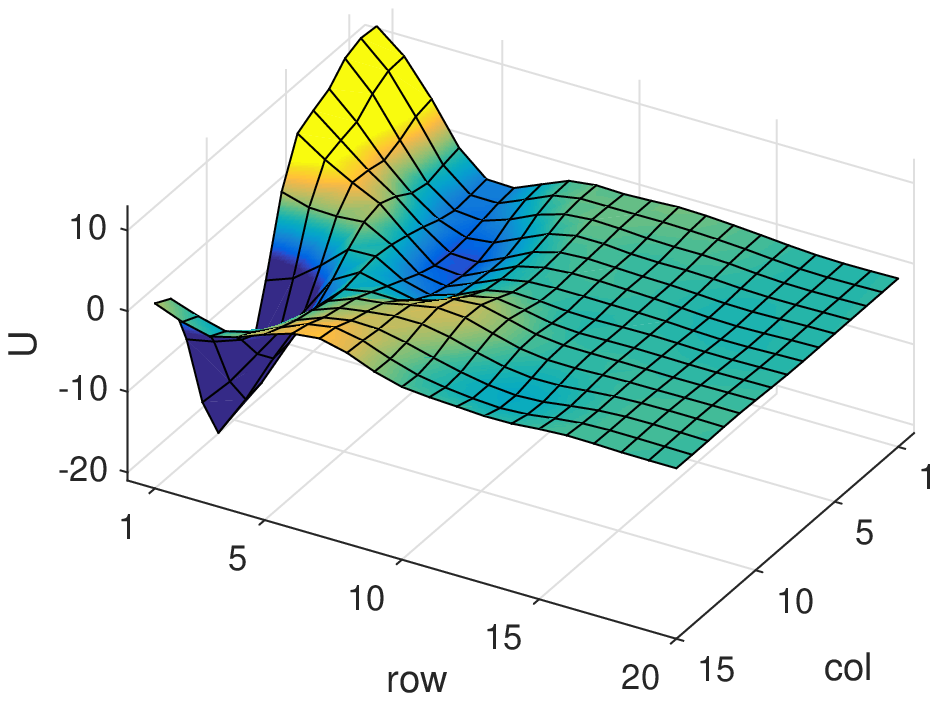}\hspace{1cm}
    \includegraphics[width=6cm]{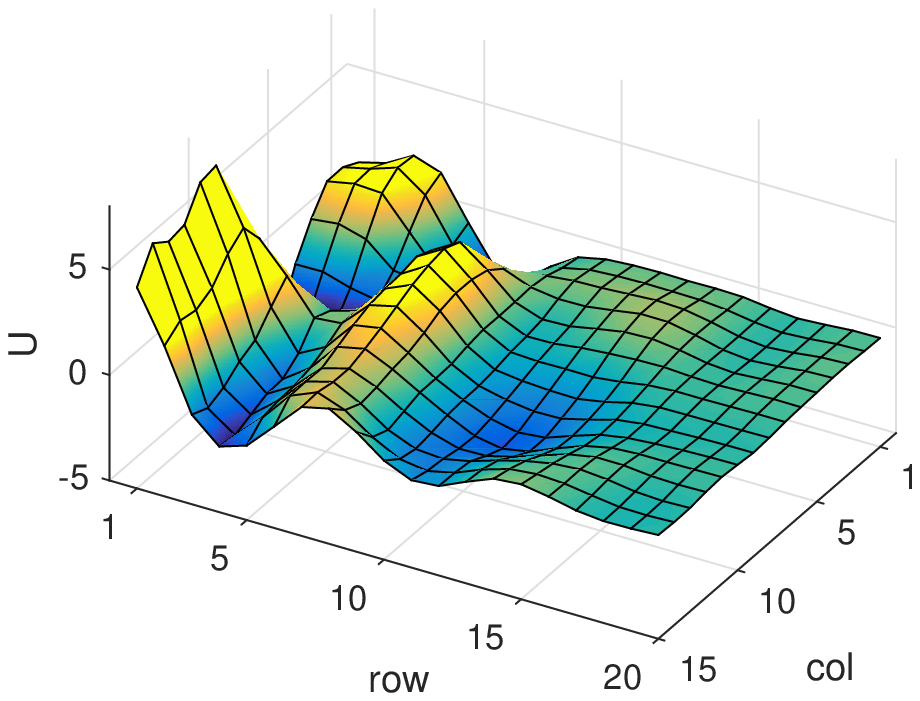}
  \end{center}
  \caption{Piecewise-linear-kernel discrete-time Urysohn operators of the two-Urysohn model shown as elements of matrix $U$. The first and the second operators are shown in the left and the right subfigures.}
  \label{fig:USurf}
\end{figure}

To see an advantage of the two-Urysohn model over the single Urysohn model, the object must be complex enough, such that the single Urysohn model is not capable of providing an accurate description. That is achieved by the shape of the kernels of generated operators. They are designed to be smooth, but with an irregular shape and with multiple local maxima and minima points, as shown in figure \ref{fig:USurf}.

For the purpose of this paper, the details regarding the operator and the input generation are omitted in the text. The full C\# source code of the example available as the supplementary information \cite{U2U}, from which the details can be obtained.

The original models and the input data are generated randomly at each execution of the test programme. The average error for $10$ experiments was $E = 3.7\% \pm 1.3\%$, where the $95\%$ confidence interval is indicated.

\subsection{Urysohn block followed by rectifier}
\label{sec:rectifier} 

The aim of this example is to show how certain information regarding the object can facilitate the identification. The considered system is a Urysohn block followed by a known nonlinearity, which is a rectifier. 

The idea of the numerical experiment is the same as in the previous example: the generation of the original model, the generation of the training and the validation datasets, the identification of the new model using the training dataset and the calculation of the error using the validation dataset. The noise is not added to the data, since the aim of the example is to estimate the accuracy of the identification method and not the error filtering capabilities. Again, for the purpose of the paper, the details are omitted in the text and can be found directly in the C\# source code of the example available as the supplementary information \cite{Rectifier}.

The nonlinearity is expressed by $z_i = \left|y_i\right|$ and cannot be uniquely inverted. There are two possible arguments $y_i$ for each known function value $z_i$. Following algorithm \ref{alg:algTwo}, intermediate output $\hat{y}$ is computed at step \ref{item:iter1} of each iteration of the identification process. It can be used to select the argument of the non-linear function (rectifier) as the closest to $\hat{y}$. In this case, at step \ref{item:iterI} of algorithm \ref{alg:algTwo}, such $\Delta y$ is taken that $z_i = \left|\hat{y}_i + \Delta y\right|$. The identification is conducted as a real-time process. 

The average error for $10$ experiments was $E = 0.7\% \pm 0.1\%$.

\subsection{Urysohn block followed by relay}
\label{sec:relay} 

The aim of this example is to show that algorithm \ref{alg:algTwo} leads to an efficient identification of objects with complex nonlinearities, where there is a significant loss of information when the intermediate variable is converted to the output. 

The idea of the numerical experiment follows the previous examples: the generation of the model and the datasets, the identification of the new model and the error calculation. The C\# source code of the example available as the supplementary information \cite{Relay}, from which the details can be obtained.

The original system is comprised of an Urysohn block, followed by a relay with two possible outputs: $+1$ or $-1$. The Urysohn block results in intermediate variable $y_i$ with both positive and negative values. The relay is the signum function, $z_i = \operatorname{sgn}\left(y_i\right)$.

After the original input-output dependence is generated, the system is identified as the canonical model with the piecewise-constant-kernel Urysohn block, equation \eqref{eq:UryshonDiscr}, and the piecewise-constant nonlinearity. Following remark \ref{rem:realTime}, the input and the output sequences are recorded and the identification is run several times through the data. The identification process switches at each run between the Urysohn operator and the nonlinearity, i.e. either step \ref{item:iter2m} or step \ref{item:iter2} is executed during one run and another step is executed during a subsequent run. 

Similarly to the example of section \ref{sec:rectifier}, at the run when the Urysohn operator is updated, $\Delta y$ varies from one iteration to another. After intermediate output $\hat{y}$ is computed at step \ref{item:iter1}, a smallest $\Delta y$, such that $z_i = \operatorname{sgn}\left(\hat{y}_i + \Delta y\right)$, is selected. At the run when the nonlinearity is updated, $\Delta y = 0$.

The average error for $10$ experiments was $E = 0.1\% \pm 0.07\%$. Identification of such system by a single Urysohn model gives a significantly more inaccurate result.


\subsection{Physical object}
\label{sec:WienerHammer} 

The benchmarking would be incomplete without applying the suggested technique to a real physical object. For this purpose, a publicly available experimental dataset has been taken \cite{GroundVibration,Schoukens2009}. The investigated object is a Wiener-Hammerstein object with static nonlinearity sandwiched between two linear blocks. The large dataset size (188K points) allows using the first half of the dataset for the training and the second half for the validation. The identification process succeeded with an error below $1\%$ on the training data, and when the model has been applied to the validation (unseen) data, it resulted in an approximate error of $1.5\%$.

It has been noticed by examining the data that statistical properties of the input signal are different at the beginning and at the end of the recorded dataset. This means that the input is not exactly a stationary random process. For such processes or objects with gradually changing states, the prediction of the output can be combined with the updating of the model in real time, following algorithm \ref{alg:algTwo}. In this case, the error on the validation data reduces to $0.7\%$. The model is updated after the error is computed; thus, the error is indeed computed for the unseen data and both operators are updated right after the error calculation, before the subsequent iteration. The full C\# source code as well as the data are available as the supplementary information \cite{WienerHammer}.

\section{Conclusions}

This paper introduces a generalisation of well-known block-oriented models --- the canonical model, given by equations \eqref{eq:canonical1}-\eqref{eq:canonical2}, and a way of identifying it. It has been shown that any deterministic SISO (single input, single output) object can be approximated by the canonical model with a predefined accuracy.

The canonical model can describe objects with any reversible (stateless) nonlinearities, which are nonlinearities that return to the original state after the input to the nonlinearity is rolled back. Therefore, the limitation of the model is that it cannot reproduce the behaviour of objects with hysteresis, friction or clearance --- non-reversible nonlinearities, as for their description a state parameter is required. 

The identification algorithm allows identifying models in real time, since both operators are updated for each new pair of input-output values using a short list of previously recorded values. Only the approximations to the identified operators and the most recent input and output values are kept in memory. 

This paper provides a number of examples with the source codes. As seen from the results, even in the case of complex nonlinearities, where there is a significant loss of information when the intermediate variable is converted to the output, it is possible to identify the model with an error below $1\%$ when tested on an independent dataset. Furthermore, an example of a real physical object has been identified with an error around $1\%$. Providing the source codes as the supplementary information, this paper offers practical tools and methods for engineers and lab researchers working in the field of dynamic systems' modelling.





\section*{Acknowledgements}

The authors of this paper would like to express a great appreciation for the data \cite{GroundVibration}, which has been collected by enthusiastic researchers, and which has been published specifically for the purpose of testing identification algorithms. 

\bibliographystyle{unsrt}
\bibliography{refs}

\end{document}